\documentclass[12pt]{article}
\usepackage{amssymb}
\usepackage{amsmath}
\usepackage[dvips]{graphicx}

\newtheorem{theorem}{Theorem}[section]

\newtheorem{proposition}[theorem]{Proposition}

\newcommand{\calH}{{\cal H}}
\newcommand{\ep}{\varepsilon}
\newcommand{\dN}{{{\bf N}}}
\newcommand{\dR}{{{\bf R}}}
\newcommand{\E}{{{\bf E}}}

\newcommand{\prob}{{{\bf P}}}

\begin{document}
\title{On the Optimal Amount of Experimentation in Sequential Decision Problems%
\thanks{We thank Ehud Lehrer for the discussions we had on the subject.}}
\author{Dinah Rosenberg\thanks{Laboratoire d'Analyse G\'{e}om\'{e}trie et
Applications, Institut Galil\'{e}e, Universit\'{e} Paris Nord,
avenue Jean-Baptiste Cl\'{e}ment, 93430 Villetaneuse, France; and
laboratoire d'Econom\'etrie de l'Ecole Polytechnique, 1, rue
Descartes 75005 Paris, France. e-mail:
dinah@zeus.math.univ-paris13.fr}\ , Eilon Solan\thanks{The School of
Mathematical Sciences, Tel Aviv University, Tel Aviv 69978, Israel.
e-mail: eilons@post.tau.ac.il} \ and Nicolas
Vieille\thanks{D\'{e}partement Finance et Economie, HEC, 1, rue de
la Lib\'{e}ration, 78 351 Jouy-en-Josas, France. e-mail:
vieille@hec.fr}}
\maketitle

\noindent
2000 Mathematics Subject Classification:
62C10, 60G99,
93E35.

\begin{abstract}
We provide a tight bound on the amount of experimentation under the optimal strategy in sequential decision problems.
We show the applicability of the result by providing a bound on the cut-off in a one-arm bandit problem.
\end{abstract}

Keywords: experimentation, sequential decision problems, optimal strategy.

\section{Introduction}

A basic issue faced by the statistician in sequential decision
problems is the trade-off between the cost of pursuing the
experimentation and the informational benefit from doing so. For
instance, in bandit problems, the decision maker chooses whether
to pull an apparently optimal arm, or to pull some seemingly
poorer one, in the hope of thereby getting valuable information.

Such problems lead to unwieldy analytical problems, rarely
amenable to closed-form solutions, which is arguably one reason
why sequential methods are still seldom relied upon in practice
(see Lai (2001), Armitage (1975)). For bandit problems, while the
optimal strategy is well characterized and consists in pulling the
arm with highest dynamic allocation index (Gittins
and Jones (1974), Gittins (1979)),  the explicit computation of these indices is
rarely feasible, except for very specific cases where the risky
arm yields a Bernoulli payoff (see for instance Bradt, Johnson and
Karlin (1956), Feldman (1962), Woodroofe (1979), Berry and
Fristedt (1985)).

Over the years, a number of approaches have been pursued: (i)
computing approximate solutions of the corresponding dynamic
programming equation, as in Berry (1972) or Fabius and van Zwet
(1970); (ii) relying on close-by problems for which explicit
solutions are known, as in Lai (1987); (iii) using extensively
numerical computations, as in Lai (1988, for sequential testing of
composite hypotheses); (iv) designing \textit{ad hoc} policies,
sometimes investigating their performance numerically, as in
Cornfield, Halperin  and Greenhouse (1969), Berry and Sobel
(1973), Berry (1978) and, more recently,   (v) finding explicit
\textit{a priori} bounds, as in Brezzi and Lai (2000).

This note contributes to the last category. Motivated by economic
applications, (see, e.g. Dixit and Pindyck (1994), Bolton and
Harris (1999), Bergemann and V\"{a}lim\"{a}ki (2000), Keller, Rady and
Cripps (2005), Rosenberg, Solan and Vieille (2007)), we consider
general Bayesian, discounted sequential problems. The parameter
$\theta$ has an initial distribution $\prob$ (the belief of the
economic agent). The agent repeatedly receives some information, chooses an action from a set
$A$, and get a possibly unobserved instantaneous reward
$u(\theta,a)$. Future gains are
discounted by a discount factor $\delta\in (0,1)$. Given a
decision rule $\sigma$, and a stage $n$, we define the
\textit{amount of experimentation} in stage $n$ to be the
difference $\Delta_n$ between the currently highest reward, and the
current reward obtained when using $\sigma$.

We show that, for every \textit{optimal} decision rule, the
expected value of $\displaystyle \sum_{n=1}^{+\infty}\Delta_n$ does
not exceed $C\delta /(1-\delta)$, where $C$ is a bound\footnote{In particular, $\displaystyle\sum\Delta_n
<\infty$ a.s., hence any optimal decision rule eventually stops to experiment.}  on the
reward function
$u$. The bound is valid irrespective of the
prior belief $\prob$, and no matter how information flows in to
the decision maker.
This result was used in Rosenberg, Solan and Vieille (2009) to show that
the limit payoff of neighbors in connected social networks coincides,
and to provide conditions that ensure concensus.

We next show, by means of an example, that this bound is tight. We
also illustrate how to use this bound in practice to derive a priori
estimates for specific sequential problems. For simplicity, we focus
on an instance of a one-arm bandit problem, for which no explicit
solution is available, and give an estimate of the optimal boundary
in the associated optimal stopping problem.
In contrast to Brezzi and Lai (2000),
who provide a bound on the Gittins' index in bandit problems,
our bound is on the cut-off of the optimal strategy.

\section{Setup and Results}

The parameter set\footnote{In spite of
the qualifier ``parameter", our decision problems are
non-parametric, since the space $\Theta$ is fully general.}
is a measurable space $(\Theta,\mathcal{A})$,
endowed with a prior distribution $\prob$.
At
each stage $n\geq 1$, a decision maker first gets an observation
drawn from a (measurable) set $S$, then chooses an action $a$ out
of a (compact metric) set $A$, and gets a reward $u(\theta,a)$.
The decision maker discounts future rewards at the rate $\delta\in
[0,1)$. The reward function $u:\Omega\times A\to \dR$ is (jointly)
measurable, and continuous w.r.t. $a$. In addition, we assume that
the highest reward $\overline{u}\colon\theta\mapsto \max_{a\in
A}u(\theta,a)$ and the lowest reward
$\underline{u}\colon\theta\mapsto \min_{a\in A}u(\theta,a)$ have
finite expectation.

We stress that we place no restriction whatsoever on the nature of observations:\footnote{Beyond
the minimal, technical assumption that
the observation in stage $n$ is drawn according to a transition
probability from $\Theta\times (S\times A)^{n-1}$ to $S$.}
e.g., they may depend, possibly in a random way, on the parameter
$\theta$, and on past observations and actions; they may or may not reveal past rewards;
and they may be independent or not.

Note that we assume that the current reward is a
\emph{deterministic} function $u(\omega,a)$ of the parameter
$\omega$ and of the action $a$. This assumption is made without
loss of generality. Statistical models such as multi-armed bandit
problems, where the decision maker observes her current reward
that randomly depends on $\theta$ (and on $a$), can be cast into
the above framework. Indeed, it suffices to re-label such a random
reward as the ``observation", and to define the reward to be the
expectation of the ``observation". Such a change does not affect
the optimal decision rules, nor the optimal value of the problem.

For a  decision \footnote{That is, a sequence $(\sigma_n)$
of measurable functions, where  $\sigma_n:(S\times A)^{n-1}\times
S\to A$ is the decision in stage $n$.} rule $\sigma$,
$\prob_\sigma$ is the joint distribution of $\theta$ and of the
infinite sequence of observations and decisions. Expectation
w.r.t.~$\prob_\sigma$ is denoted by $\E_\sigma$.

We focus on the amount of experimentation that optimal decisions
entail. To be specific, let a decision rule $\sigma$ be given.
Given a stage $n$, we denote by $\calH_n$ the information
available at stage $n$, that is, the $\sigma$-field induced by
past observations and actions. When using the decision rule
$\sigma$ prior to stage $n$, the expectation
$\E_\sigma[u(\theta,a)|\calH_n]$ is the expected reward when
choosing $a$ in stage $n$, given all available information, and
$\overline{u}_n:=\max_{a\in A}\E_\sigma[u(\theta,a)|\calH_n]$ is
the myopically optimal reward. Thus, letting $a_n$ denote the
action of the decision maker in stage $n$,
$u_n=\E_\sigma[u(\theta,a_n)|\calH_n]$ is the actual reward that
the decision maker expects to get in stage $n$, when following
$\sigma$. The difference
$\Delta_n:=\overline{u}_n-u_n$
provides a measure of the degree of
experimentation  performed in stage $n$. The infinite sum
$\displaystyle \sum_{n\geq 1}\Delta_n$ therefore measures the overall
amount of experimentation.

\begin{theorem}\label{theorem rate-convergence}
For any optimal decision rule $\sigma$, one has

\[\E_{\sigma}\left[\sum_{n\geq 1}\Delta_n\right]\leq \left(\E\left[\overline{u}\right] - \E\left[\underline{u}\right]\right)\times \displaystyle
\frac{\delta}{(1-\delta)}.\]
\end{theorem}

Beyond quantitative implications, this bound also yields
qualitative implications. Consider for instance a multi-arm bandit
problem. For simplicity, assume that the types of the various arms
are first drawn, and that each arm then yields a sequence of
rewards, which is conditionally i.i.d. given its type. For
concreteness, assume that with probability 1 over the types, the
expected outputs of the arms are all distinct.

Observe that, whenever the decision maker pulls a specific arm
infinitely often, she eventually learns the type of this arm.
Therefore, whenever the decision maker pulls \textit{two} specific arms
infinitely often, she eventually learns both types. Since one of
these two arms is ``better" than the other, this implies that the
sequence $(\Delta_n)_{n\geq 1}$ then does not converge to zero. By Footnote
1, this event must have probability 0, for every optimal decision
rule. In other words: any optimal allocation rule samples finitely
often all arms but one.
This provides an alternative proof of Theorem 2 in Brezzi and Lai
(2000).\footnote{Brezzi and Lai (2000) assumes
that the states of the different arms are independent. Our argument dispenses with this assumption.}

\vskip\baselineskip We next show that the bound in Theorem
\ref{theorem rate-convergence} is tight.

\begin{proposition}\label{first example}
For every $\ep$ and for every discount  factor $\delta$, there is a
decision problem with an optimal decision rule  $\sigma$ such that
$\E_{\sigma}[\sum_{n\geq 1}\Delta_n]\geq
(\E[\bar{u}]-\E\left[\underline{u}\right])\displaystyle\times
\frac{\delta}{(1-\delta)}\times \left(1-\ep\right)$.
\end{proposition}

The decision problem in Proposition \ref{first example} depends both
on $\ep$ and on $\delta$. The next proposition improves in this
respect, at a slight cost in the speed of convergence.
 In this statement, and given $\ep>0$, we denote by $N(\ep)$ the (random) number of stages
in which $\Delta_n$ is at least $\ep$: $N(\ep):=|\{n\geq 1: \Delta_n\geq
\ep\}|$. Plainly, $\displaystyle \sum_{n\geq 1}\Delta_n\geq \ep
N(\ep)$ for every $\ep > 0$.

\begin{proposition}\label{second example}
There is a decision problem such that for every $\delta> 2/3$
there is a unique optimal decision rule $\sigma$ that satisfies
\[\lim_{\ep\to 0}\ep^\alpha \E_{\sigma}\left[N(\ep)\right]=+\infty,\mbox{ for every }\alpha <1.\]
\end{proposition}
That is, as $\ep$ decreases, the expected number
$\E_{\sigma}[N(\ep)]$ of experimentation stages
 increases faster than $1/\ep^\alpha$, for every $\alpha <1$.

\section{Proofs}
\label{proofs}

\subsection{Proof of Theorem \ref{theorem rate-convergence} }
\label{proof theorem rate}

Consider an optimal decision rule $\sigma$. Set $\displaystyle
Y_n:=(1-\delta)\sum_{k=n}^{+\infty}\delta^{k-n}\E_{\sigma}[u_k\mid
\calH_n]$: $Y_n$ can be interpreted as the \textit{continuation}
reward  under the optimal decision rule  (discounted back to stage
$n$). Since $u_k \leq \E_{\sigma}[\overline{u}\mid \calH_k]$ for
all $k\geq n$, one has $\E_{\sigma}[Y_n] \leq \E[\overline{u}]$.

Since one option available to the decision maker, from stage $n$
on, is to ignore all future observations, and to keep choosing the
action that was myopically optimal in stage $n$, we have
\begin{equation}
\label{equ 2.1} Y_n\geq \overline{u}_n.
\end{equation}
Now, rewrite $Y_n$ as
\begin{eqnarray}
Y_n &=& (1-\delta)u_n+\delta \E_{\sigma}[Y_{n+1}\mid \calH_n]\nonumber\\
\label{basic} &=& (1-\delta) \left(\overline{u}_n-\Delta_n \right) + \delta
\E_{\sigma}[Y_{n+1}\mid \calH_n].
\end{eqnarray}
From (\ref{equ 2.1}) and (\ref{basic}) we obtain:
\[ \overline{u}_n \leq (1-\delta) \left(\overline{u}_n-\Delta_n\right) + \delta \E_{\sigma}[Y_{n+1}\mid \calH_n], \]
so that after cancelling $\overline{u}_n$ from both sides and dividing by $\delta$,
\begin{equation}\label{basic2}
\overline{u}_n\leq
\E_{\sigma}\left[Y_{n+1}\mid\calH_n\right]-\frac{\Delta_n(1-\delta)}{\delta}.
\end{equation}
Substituting (\ref{basic2}) into (\ref{basic}), we obtain
\begin{eqnarray*}
Y_n &\leq&
(1-\delta)\left(\E_{\sigma}\left[Y_{n+1}\mid\calH_n\right] -
\Delta_n\left(\frac{1-\delta}{\delta}+1\right) \right)
 + \delta \E_{\sigma}\left[Y_{n+1}\mid\calH_n\right]\\
&\leq& \E_{\sigma}\left[Y_{n+1}\mid\calH_n\right]  -
\frac{1-\delta}{\delta}\Delta_n .
\end{eqnarray*}
Taking expectations, summing over $n=1,\ldots,k$, using
$\E\left[\underline{u}\right] \leq \E_{\sigma}[Y_n] \leq \E\left[\overline{u}\right]$, and taking the limit
as $k$ goes to infinity, we obtain
\[ \E_{\sigma}\left[\sum_{n\geq 1}\Delta_n\right] \leq \left(\E[\overline{u}]-\E[\underline{u}]\right) \times \frac{\delta}{(1-\delta)}, \]
as desired.

\subsection{Proof of Proposition \ref{first example}}

Fix $\delta > 0$.
Note that if the statement holds for $\ep_0$, then it holds for every $\ep > \ep_0$.
We will prove that the statement holds for $\ep = 1/m$, for any natural number $m > 1/\delta$.
Let $\Theta = \{\theta_1,\theta_2,\ldots,\theta_m\}$ and $A =
\{a_0,a_1,\ldots,a_m\}$ contain $m$ and $m+1$ elements respectively.
The prior belief on $\Theta$ is uniform, and
 the reward function is given by :
\begin{eqnarray}
u(\theta_k,a_k) &=& 1, \ \ \ \ \ k=1,\ldots,m,\\
u(\theta_k,a_l) &=& 0, \ \ \ \ \ k=1,\ldots,m,l \neq k,\\
u(\theta_k,a_0) &=& 0, \ \ \ \ \ k=1,\ldots,m.
\end{eqnarray}
Thus, once the parameter is inferred with certainty, there is a
unique optimal action, whereas \textit{ex ante}, $a_1,\ldots,a_m$
are all myopically optimal, while $a_0$ is $(1/m)$-suboptimal.

Information is provided to the decision maker according to the
following rules:
if the decision maker has chosen $a_0$ in all previous stages,
the true parameter is revealed with probability $c
:=\frac{(1-\delta)}{\delta(m-1)}<1 $;
if the decision maker did not choose $a_0$ in all previous stages, no information is revealed,
that is, no observation is made.
Suppose the decision maker chooses $a_0$ until the state of the
world is revealed, and then switches to the optimal action. The
expected reward $A$ satisfies $A = c\delta + (1-c)\delta A$, so that $A = \frac{c\delta}{1-(1-c)\delta}$.
Substituting $c =\frac{(1-\delta)}{\delta(m-1)} $ we
obtain that the expected reward is $1/m$, so that this strategy is
optimal. However, for $\ep = 1/m$ one has:
\[ \E_\sigma\left[\sum_{n\geq 1} \Delta_n\right] = \E_\sigma[\ep N(\ep)] = \frac{\ep}{c} = \frac{m-1}{m}\frac{\delta}{1-\delta}.\]
Since $\bar{u}=1$ and
$\underline{u}=0$ we  get the desired result.

\subsection{Proof of Proposition \ref{second example}}

We provide an example within the class of Gaussian models. Set
$\Theta=\dR$, and let the action set $A=\dR\cup\{-\infty,+\infty\}$
be the set of extended real numbers, endowed with the usual
topology. The reward function $u(\theta,a)$ is equal to one if $a\in
\dR$ and $|\theta-a|\leq 1$, and equal to zero otherwise.

Given a normal distribution $\mu$ with precision $\rho$ (that is, with variance
$1/\rho$), define
$\bar{u}(\rho)$ to be the highest reward that the decision maker
may achieve, when holding the belief $\mu$. Observe that
$\bar{u}(\rho)$ does not depend on the mean of $\mu$. Plainly, the
map $\rho\mapsto \bar{u}(\rho)$ is continuous and increasing, with
$\lim_{\rho\to 0}\bar{u}(\rho)=0$, and $\lim_{\rho\to
+\infty}\bar{u}(\rho)=1$.

The signalling structure of the decision problem is designed in
such a way that the decision maker's belief is always a normal
distribution. In addition, she keeps receiving additional
information about $\theta$ as long as she follows a pre-specified
sequence of suboptimal actions.

To be specific, let $(\ep_n)_{n \geq 1}$ be a decreasing sequence
of positive numbers that satisfies (i) $\sum_{n=1}^\infty \ep_n
\in (1/2,1)$, (ii) $\ep_n n^\beta \to +\infty$, for every
$\beta>1$, and\footnote{For instance, choose
$\ep_n=\frac{(n\ln^2n)^{-1}}{\sum_{k=1}^\infty (k\ln^2k)^{-1}}$
for $n$ sufficiently large.}
(iii) $\frac{\ep_{n-1}}{\ep_n} >
\frac{2}{3}$.
The sequence $(\rho_n)_{n\geq 1}$ is defined recursively by the condition
\[\bar{u}(\rho_1+\cdots+\rho_n)=\ep_1+\cdots+\ep_n.\]

Let the prior distribution $\prob$ be a normal distribution with
precision $\rho_1$, and let $(\xi_n)_{n\geq 2}$ be a sequence of
independent normally distributed variables with precision $\rho_n$, and
independent from $\theta$.

Observe that, in the absence of any information about $\theta$,
the decision maker's myopically optimal reward is
$\bar{u}(\rho_1)=\ep_1$. We set $a_1=+\infty$. On the other hand,
if she receives the observations $s_k:=\theta+\xi_k$, $k=2,\cdots,
n$ ($n\geq 2$), her belief over $\theta$ is normally distributed,
with precision $\rho_1+\cdots+\rho_n$. Hence, her myopically
optimal reward is
$\bar{u}(\rho_1+\cdots+\rho_n)=\ep_1+\cdots+\ep_n$, and there is
an action $a_n$ (which depends on $s_2,\ldots,s_n$), which yields
an expected reward equal to $\ep_1+\cdots,+\ep_{n-1}$.

We now define the information received by the decision maker:
\begin{itemize}
\item   Prior to stage 1, the decision maker receives no observation;
\item Prior to stage 2, she receives the observation $s_2=\theta+\xi_2$
if she played $a_1=+\infty$ at the first stage, and no observation
otherwise; \item Prior to stage $n>2$, she receives the
observation $s_n=\theta+\xi_n$ if she played $a_1,a_2,\ldots,
a_{n-1}$ at the previous stages. Otherwise, she receives no
observation.
\end{itemize}

Playing the sequence $(a_n)$ of actions is the unique optimal
decision rule. Indeed, if the decision maker first deviates from
that sequence at stage $k\geq 1$, she receives no further
information, hence her optimal reward in all later stages is
$\ep_1+\cdots+\ep_k$; if she sticks to the
sequence $(a_n)$, her continuation reward (discounted back to
stage $k$) is
\[ (1-\delta) \sum_{n=k}^\infty \delta^{n-k} (\ep_1+\cdots+\ep_{n-1}). \]
By (iii), this reward is higher than $\ep_1+\cdots+\ep_k$.

Note that $a_n$ is (myopically) $\ep_n$-optimal, for each $n\geq
1$. Since the sequence $(\ep_n)$ is decreasing, there are exactly
$n$ rounds in which the decision maker does not play a myopically
$\ep_n$-optimal action, so that by (ii) $(\ep_n)^\alpha N(\ep_n) =
n(\ep_n)^\alpha$ converges to infinity for every $\alpha < 1$.

\section{Application}

We here illustrate how Theorem \ref{theorem rate-convergence} can
be used to derive \textit{a priori} bounds on the optimal decision
rules in specific decision problems. Since our goal is here purely
illustrative, we restrict ourselves to the analysis of a specific
one-arm bandit problem, where the risky arm has two possible
types, a good type and a bad type, and observations are i.i.d. In such a problem, the optimal
decision rule consists of pulling the risky arm as long as the
posterior probability assigned to the good type exceeds a specific
cut-off, and then in switching permanently to the safe arm.

We set the problem so as to depart as little as possible from a
Bernoulli problem, for which a closed form expression for the
optimal cut-off is known. We also make no attempt at optimizing
our final bound.

The type $\theta$ of the risky arm takes values in the two-point
set $\{\theta_0,\theta_1\}$. Both types are \textit{ex ante}
equally likely. The safe arm yields zero.
Given $\theta=\theta_i$, the risky arm may yield three different
rewards, $a,b$ and $c$, with probabilities $p^i_a,p^i_b$ and
$p^i_c$.  These probabilities are such that (i) the expected
reward of the risky arm is $1$ if $\theta=\theta_1$, and $-1$ if
$\theta=\theta_0$; (ii) one has $\displaystyle
\ln\frac{p^1_a}{p^0_a}=\alpha$, $\displaystyle
\ln\frac{p^1_b}{p^0_b}=2\alpha$, and $\displaystyle
\ln\frac{p^1_c}{p^0_c}=-\alpha$, for some $\alpha>0$.

Denote by $\pi_n$ the posterior belief that $\theta=\theta_1$,
based on all observations prior to stage $n$, and let
$Z_n=\displaystyle\ln\frac{\pi_n}{1-\pi_n}$ be the log-likelihood
ratio. Conditional on $\theta=\theta_0$, the sequence $(Z_n)$
follows a random walk, which moves up by  $\alpha$ (with
probability $p^0_a$), by $2\alpha$, or moves down by $\alpha$
between any two stages.

The optimal decision rule consists in pulling the risky arm until
the first stage $\sigma^\ast$ where $Z_n=-k^\ast\alpha$, for some
$k^\ast\in \dN$, and then in pulling repeatedly the safe arm. We
will derive an upper bound on $k^\ast$ using Theorem \ref{theorem
rate-convergence}.

The amount of experimentation in stage $n$ is
$\Delta_n=\max\{0,1/2-\pi_n\}$. For $k<k^\ast$, let $N(k)$ be the
number of passage of the sequence $(Z_n)$ at the level $-k\alpha$,
and denote by
$\ep(k)=1/2-\displaystyle\frac{e^{-k\alpha}}{1+e^{-k\alpha}}$ the
corresponding value of $\Delta_n$. Thus,
\begin{equation}\label{identity}\sum_{n=1}^{+\infty}\Delta_n=\sum_{k<k^\ast}\ep(k)N(k).\end{equation}

Observe now that whenever $Z_n=-k\alpha$, the expected number of
visits (including stage $n$) to $-k\alpha$ before $Z_n$ moves
below $-k\alpha$ is $1/(1-p^0_a)$. On the other hand, it is then
the case that the sequence $(Z_n)$ moves down to $-(k+1)\alpha$.
Hence, the probability that $(Z_n)$ will move back to $-k\alpha$
before hitting $-k^*\alpha$ is\footnote{This bound is
admittedly very crude.} at least $p^0_a$. Therefore,
\begin{equation}\label{identity2}\E_{\theta_0}[N(k)]\geq \frac{p^0_a}{1-p^0_a}.\end{equation}
By Theorem \ref{theorem rate-convergence} one has
$\displaystyle
\frac{1}{2}\E_{\theta_0}\left[\sum_{n=1}^{+\infty}\Delta_n\right]+\frac{1}{2}\E_{\theta_0}\left[\sum_{n=1}^{+\infty}\Delta_n\right]\leq
\frac{2\delta}{1-\delta}$. Therefore,  (\ref{identity}) and
(\ref{identity2}) yield
\begin{equation}\label{bound3}\sum_{k=0}^{k^*-1}\frac{1}{2}\frac{1-e^{-k\alpha}}{1+e^{-k\alpha}}=\sum_{k=0}^{k^*-1}\ep(k)\leq
4\frac{1-p^0_a}{p^0_a (1-\delta)},
\end{equation}
By monotonicity, the left-hand side of (\ref{bound3}) is at least
equal to
\[
\frac{1}{2}\int_0^{k^*-1}\tanh\frac{x\alpha}{2}dx
=\frac{1}{\alpha}\ln\cosh \frac{\alpha(k^*-1)}{2} \geq
\frac{1}{\alpha}\ln\frac{e^{\frac{\alpha(k^*-1)}{2}}}{2}
=\frac{(k^*-1)}{2}-\frac{\ln 2}{\alpha}.
\]
Thus,
\[ k^*\leq 4\left(1+2\frac{\ln 2}{\alpha}+2\frac{1-p^0_a}{p^0_a
(1-\delta)}\right).\]

\end{document}